\documentclass[review]{elsarticle}

\usepackage{lineno,hyperref}
\modulolinenumbers[5]

\usepackage{stmaryrd}
\usepackage{amssymb}
\usepackage{color, amsmath,amssymb, amsfonts, amstext,amsthm, latexsym}

\usepackage{amssymb, epsfig, amssymb, latexsym}
\usepackage{amsmath}
\usepackage{graphicx}
\usepackage{longtable}
\usepackage{graphicx}
\usepackage{subfigure}

\renewcommand{\phi}{\varphi}

\begin{document}

\begin{frontmatter}

\title{Detecting the maximum likelihood transition path from data of stochastic dynamic systems\tnoteref{mytitlenote}}

%\author{Min Dai\fnref{myfootnote}}
%\address{School of Mathematics and Statistics, \\Center for Mathematical Science, \\ Huazhong University of Science and Technology, Wuhan, 430074, China}
%\ead{mindai@ust.edu.cn.}
%\fntext[myfootnote]{mindai@ust.edu.cn}

\author[mymainaddress]{Min Dai}
\ead{mindai@hust.edu.cn}

\author[mysecondaryaddress]{Ting Gao\corref{mycorrespondingauthor}}
\cortext[mycorrespondingauthor]{Corresponding author}
\ead{tinggao0716@gmail.com}

\author[mymainaddress]{Yubin Lu}

\author[mymainaddress,mythirdaddress]{Yayun Zheng}

\author[myfourthaddress]{Jinqiao Duan}
\ead{duan@iit.edu}

\address[mymainaddress]{School of Mathematics and Statistics, \& Center for Mathematical Science, \\ Huazhong University of Science and Technology, Wuhan, 430074, China.}
\address[mysecondaryaddress]{Twitter, 1335 Market St \#900, San Francisco, CA 94103, USA.}
\address[mythirdaddress]{Wuhan National Laboratory for Optoelectronics,\\ Huazhong University of Science and Technology, Wuhan, 430074,China.}
\address[myfourthaddress]{Department of Applied Mathematics, College of Computing,\\ Illinois Institute of Technology, Chicago, IL 60616, USA.}

\begin{abstract}
   In recent years, the discovery of complex dynamic systems in various fields through data-driven methods has attracted widespread attention. This method has played the role of data and has become an advantageous tool for us to study complex phenomena. In this work, we propose a framework for detecting the dynamic behavior, such as the maximum likelihood transition path, of stochastic dynamic systems from data. For the stochastic dynamic system, we need to use the Kramers-Moyal formula to convert it into a deterministic form for processing, then use the extended SINDy method to obtain the parameters of stochastic dynamic systems, and finally calculate the maximum likelihood transition path. We give two examples of stochastic dynamical systems driven by additive and multiplicative Gaussian noise, and demonstrate the validity of the method by reproducing the known dynamical system behavior.
\end{abstract}


\end{frontmatter}

\linenumbers

\section{Introduction}
Stochastic differential equations are widely used to describe random phenomena in disciplines such as physics, biology, chemistry, and geophysics. For such a stochastic problem, we usually build an appropriate mathematical model based on the basic laws, and then analyze or simulate the model to obtain a characterization of the nonlinear phenomena of the problem. However, for some phenomena that are too complicated, we lack sufficient understand for them, and it is difficult to establish complete mathematical models, or the models corresponding to some phenomena are too complicated to analyze. Fortunately, with the improvement of observation technology and computing power, although there is not enough understanding of the problem, there are still many valuable observation or simulation data that can be used. Therefore, it is necessary to directly discover the dynamic system indicators from data and obtain a characterization of stochastic phenomena.

The analysis of dynamical complex behavior based on data has received extensive attention in recent years. Many authors have come up with insightful methods based on areas they are familiar with. For example, the sparse identification of nonlinear dynamics (SINDy) was used by Brunton, Kutz and et al. \cite{Bru,Rudy,Alla}, to discover the governing equation from data. Furthermore, SINDy was applied in the learning of biological networks, which effectively dealt with the problem of rational functions \cite{Man}. Zhang and Lin \cite{Zhang} used threshold sparse Bayesian regression to discover the governing physical laws from data. Moreover, there were many other methods, such as stochastic parametrization \cite{Cho}, learning informed observation geometries \cite{Yai}, the Koopman Operator \cite{Wil,Sch}, Gaussian processes \cite{Arc,Rut,Batz}, extended the SINDy methods \cite{Bon} and so on. However, most of authors focused on system identification \cite{Gar,Sch} and lacked sufficient attention to extract dynamic system indicators \cite{Wu} from data.

The maximum likelihood transition path is a significant indicator for describing the behavior of stochastic dynamic systems \cite{Zheng,Duan,Gao}, which provides information about the transition phenomena under the interaction of nonlinearity and uncertainty. We know that it is difficult to extract the accurate the maximum likelihood transition path from data under the influence of nonlinearity and uncertainty. The purpose of this paper is to devise a method to extract the maximum likelihood transition path of the stochastic dynamical system from data. We first have used the Kramers-Moyal formula \cite{Schuss} to convert the stochastic problem into the deterministic problem, then estimated the coefficients of stochastic differential equations \cite{Oksendal} by extending the SINDy method \cite{Bon} from data to calculate the maximum likelihood transition path.

The work is arranged as follows. In section 2, we introduce the theory of the article, including the Kramers-Moyal formula, extending the SINDy method and the maximum likelihood transition path. In section 3, the results of the numerical experiment by stochastic differential equations with additive and multiplicative noise are presented. In section 4, we give some conclusions and discussions.

\section{Theory}
With the advancement of technology, the combination of stochastic differential equations and data has become a powerful tool for our research. In this work, we consider the stochastic differential equation in $\mathbb{R}$ as follows
\begin{equation}\label{eq:1}
dX_t=f(X_t)dt+\sigma(X_t){dW_t},~~~  X_0=x,
\end{equation}
where $f,\sigma$ are called the drift and the diffusion of the process in $\mathbb{R}$, respectively. The vector $X_t\in{\mathbb{R}}$ denotes the state of a system at time $t$ and $W_t$ is a standard Brownian motion. For the stochastic differential equation \eqref{eq:1}, the drift $f$ and the diffusion $\sigma$ are two important parts. And fortunately, the Kramers-Moyal formula \cite{Schuss} provides us with a way to detect them from data
\begin{equation}\label{eq:2}
{f(x)}=\lim_{\Delta{t}\rightarrow{0}}{\mathbb{E}{\Bigg[\frac{(X_{{\Delta{t}}}-X_0)}{\Delta{t}}|X_0=x\Bigg]}},
\end{equation}
\begin{equation}\label{eq:3}
{\sigma^2(x)}=\lim_{\Delta{t}\rightarrow{0}}{\mathbb{E}{\Bigg[\frac{(X_{{\Delta{t}}}-X_0)^2}{\Delta{t}}|X_0=x\Bigg]}}.
\end{equation}
The above equations calculate the expectation and variance of the process $X(t)$ under the condition of position $x$ at time $t$ respectively. If the expectation and variance can be calculated from data, we can easily obtain the parameters of the stochastic differential equation.

\subsection{\textbf{{Data-Driven Parameterization of Stochastic Differential Equation} }}
We now apply the extend the SINDy algorithm \cite{Bon} by considering conditional expectation \eqref{eq:2} and \eqref{eq:3} with model $Y={\Theta\Xi}$, where $Y^1=\lim_{\Delta{t}\rightarrow{0}}{\mathbb{E}{[\frac{(X_{{\Delta{t}}}-X_0)}{\Delta{t}}|X_0=x]}}$ or $Y^2=\lim_{\Delta{t}\rightarrow{0}}{\mathbb{E}{[\frac{(X_{{\Delta{t}}}-X_0)^2}{\Delta{t}}|X_0=x]}}$. For example, suppose that the basis function is polynomial $\{1,X,X^2,X^3,\cdots\}$, then using the basis function
\begin{equation*}
\Theta(X)=
\left[
  \begin{array}{ccccc}
    1 & x_1 & x_1^2 & x_1^3 & \cdots\\
    1 & x_2 & x_2^2 & x_2^3 & \cdots\\
    \vdots & \vdots & \vdots & \vdots & \vdots \\
    1 & x_N & x_N^2 & x_N^3 & \cdots\\
  \end{array}
\right] ,
\end{equation*}
we can change Eqs. \eqref{eq:2} and \eqref{eq:3} into the following form
\begin{equation}\label{eq:4}
\begin{split}
Y^1=\lim_{\Delta{t}\rightarrow{0}}{\mathbb{E}{\Bigg[\frac{(X_{{\Delta{t}}}-X_0)}{\Delta{t}}|X_0=x\Bigg]}}={\Theta(X)\Xi^1},\\
Y^2=\lim_{\Delta{t}\rightarrow{0}}{\mathbb{E}{\Bigg[\frac{(X_{{\Delta{t}}}-X_0)^2}{\Delta{t}}|X_0=x\Bigg]}}={\Theta(X)\Xi^2},
\end{split}
\end{equation}
where $\Xi^i=[\xi_1^i,\xi_2^i,\cdots,\xi_N^i]^T,~i=1,2$ are the coefficients of the sparse vectors.

Vectors $\Xi^1$ and $\Xi^2$ are sparse vector of coefficients. They can be obtained by the following formula
$$\Xi^i=\Theta(X)^{\dag}Y^i,~~~i=1,2. $$
Brunton et al.\cite{Bru} considered such a learned model, however, their only focused on the deterministic dynamical system. In our work, we consider the stochastic dynamical system. Due to the influence of random term, the original optimization method can not meet our requirements, so we introduce a new optimization form \cite{Bon}.

\subsubsection*{\rm{\textbf{1. Binning}}}
The effect of random noise is a difficult problem in dealing with stochastic problems, so we need to preprocess the data to reduce noise. The coordinate binning is an effective method to reduce the influence of noise. To solve Eq.\eqref{eq:4}, we need to use the data $X$ from the stochastic differential equation \eqref{eq:1}. Therefore, we first place the time series data $X$ into $G$ bins,
$$\{X(t_i)\}_{i=1,\cdots,N}\mapsto\{\overline{x}_j\}_{j=1,\cdots,G},$$
with $\overline{x}_j$ representing the center of the $j$th bin. Next, we have
$$\mathbb{X}\in{\mathbb{R}^{N\times{M}}}\mapsto{\mathbb{X}_G\in{\mathbb{R}^{G\times{M}}}},$$
$$\mathbb{Y}\in{\mathbb{R}^{N}}\mapsto{\mathbb{Y}_G\in{\mathbb{R}^{G}}},$$
with the element of ${Y}_G$ representing the center of each bin.

\subsubsection*{\rm{\textbf{2. Stepwise Sparse Regressor}}}
Inspired by Brunton et al.\cite{Bru}, we will use the threshold approach to find the sparse coefficients. We first use a standard linear regression which is unconstrained to compute Eq.\eqref{eq:4} to obtain a non-sparse solution $c$. Then, let the element in the coefficient $c$ that is less than the preset threshold value $\lambda$ be $0$. Next, the remaining coefficients are continued to repeat the linear regression of the first step until no element in the coefficient $c$ is found to be less than the threshold. However, a significant drawback to this approach is need to adjust the threshold $\lambda$ appropriately. Therefore, we introduce the Cross Validation instead of threshold $\lambda$ to avoid this problem. This algorithm is called \emph{Stepwise Sparse Regressor}$(SSR)$, and the following is a specific description£º

Step $1$: The initial solution $\widetilde{\textbf{b}}$ is obtained by solving the following least squares regression equation, which may not be a sparse, $$\widetilde{\textbf{b}}=\arg\min_{\textbf{b}\in{\mathbb{R}^M}}{\|\mathbb{Y}-\mathbb{X}\textbf{b}\|_2^2}.$$

Step $2$: Taking the absolute value of all coefficients of the initial solution $\widetilde{\textbf{b}}$ and making the smallest one be $0$, $$\widetilde{b}_i=0,~~~ i=\min_m{|\widetilde{b}_m|}.$$
This approach is similar to the threshold method, and the advantage is that the number of iterations required to run is equal to the sparsity of the solution, allowing us to program accurately.

Step $3$: Least square regression calculation for the rest again $$\mathbb{Y}=\mathbb{X}[:,\widehat{j}]\widetilde{b}[\widehat{j}],$$
with $\widehat{j}$ denoting all sets except $j$ and $\mathbb{X}[:,\widehat{j}]$ representing the remaining matrix except column $j$.

Step $4$:  Repeating the above steps continuously until the optimal solution $\widetilde{c}$ indicated by Cross Validation is reached.

Note that we simply write SSR as $$\emph{SSR}(\mathbb{X},\mathbb{Y})_q,$$ indicating that running algorithm $q$ times can obtain the solution $\textbf{b}$. This means that the solution has $q$ zeros, called $q$-sparse.

\subsubsection*{\rm{\textbf{3. Cross Validation}}}
Cross Validation (CV) is a traditional statistical verification technique that is applied in interdisciplinary fields such as model selection and hyperparameter selection. In this work, we find the optimal solution through cross validation techniques, which has the advantage that the number of iterations required for Stepwise Sparse Regression can be determined.

There are many parametric models for a given data set $\mathcal{D}$, we however want to find the most suitable mathematical model. For Cross Validation, the data set is usually divided into training set and test set. Then we obtain the cross-validation (CV) score by computing the deviation $\delta$ of the test set on the training model, which measures how to balance accuracy and predictability in the model. We finally select the parameter model corresponding to the minimum value of the CV score as the optimal model.

In addition, we know that each iteration of the algorithm will produce the solution with different degrees of sparsity $$\{SSR_q\}_{q=1,\cdots,M},$$ we then can obtain different models, and calculate the CV scores $\delta[SSR_q]$.

There are many types of cross validation methods. Here, we will partition the data set into $k$ parts and then use each part as a test set in turn to train the parameter model. Firstly, we randomly divide the data set $\mathcal{D}$ into $k$ parts with the same size, that do not intersect each other, i.e., $\bigcup_i{B_i}=\mathcal{D}$, $B_i\cap{B_j}=\emptyset$. Moreover,$$\mathbb{X}_{B_i}=\mathbb{X}[p_{B_i},:], ~~~p_{B_i}=\bigcup_{p\in{B_i}}p,$$
with $\mathbb{X}[p_{B_i},:]$ denoting the $p_{B_i}$th row of the matrix. We can next define the CV score as follows
\begin{equation}\label{eq:5}
\delta^2[SSR_q]=\frac{1}{k}{\sum_{i=1}^{k}}{\|\mathbb{Y}_{B_i}-\mathbb{X}_{B_i}\cdot{SSR(\mathbb{X}_{C_i},\mathbb{Y}_{C_i})_q}}\|_2^2,
\end{equation}
\begin{equation}\label{eq:6}
C_p=\bigcup_{i\neq{p}}{B_i},
\end{equation}
with $SSR(\mathbb{X}_{C_i},\mathbb{Y}_{C_i})_q$ representing the parameter coefficients which are $q$-sparse obtained by running $SSR$ algorithm on training set $C_i$. The cross validation score has a good measure of the accuracy, predictability and sparseness of the model. We hope that the model results corresponding to low cross validation score with variable sparsity are better. The drift and diffusion terms can be learned in the above way, thus we can extract stochastic model \eqref{eq:1} and some dynamic system behaviors.

We use the extended SINDy method in this work. On the one hand, the coordinate binning effectively reduces the impact of noise. In addition, we cleverly use the CV method to avoid the problem of adjusting the appropriate threshold, and find the number of iterations precisely, which improves the accuracy of the results.

\subsection{\textbf{The maximum likelihood transition path}}
A good tool for describing the dynamic behavior of stochastic dynamic systems is the maximum likelihood transition path \cite{Zheng}. Therefore, in this work, we will use it to characterize the information we extract from the data and determine whether it is accurate or not. The maximum likelihood transition path is a manifestation of transition of a dynamic system from one state to another. In other words, for $t\in[0,t_f]$, $x,x_0,x_{f}\in\mathbb{R}$, when conditions $X(0)=x_0$ and $X(t_f)=x_{f}$ are given, we then suppose that the function $\mathcal{P}_{A}(x,t)$ which is conditional probability density exists (where $A$ indicates these two-point conditions) and expresses as follows:
\begin{equation}\label{eq:7}
\begin{split}
\mathcal{P}_{A}(x,t)&=p(X(t)=x|X(0)=x_0;X(t_f)=x_{f})\\
&=\frac{Q(x_{f},t_f|x, t)Q(x,t|x_0,0)}{Q(x_{f},t_f|x_0,0)},
\end{split}
\end{equation}
with transition probability density $Q$, which is the solution of Fokker-Planck equation.

Next, we will give a brief introduction to $Q$. First of all, we will introduced two expressions of probability density. One is $p(X(t)=u)$ ,which we use to express the probability density of SDE solution $X(t)$ at $X(t)=u$ in $\mathbb{R}$, the other is what we call the transition density $Q(u;t|\xi,s)$, which is defined on $\mathbb{R}\times[0,t_f]\times\mathbb{R}\times[0,t_f]$.  Furthermore, there is a relationship between $p$ and $Q$. For instance, for $0\leqslant s<t \leqslant t_f$, given $X(s)=v$, then the density of $X(t)$ can be expressed by $Q(u;t|\xi;s)$ at $X(t)=u$. That is to say, we have the following expression
\begin{equation*}
Q(u;t|\xi;s)=p(x(t)=u|x(s)=v).
\end{equation*}
For every $x_0\in\mathbb{R}$, we assume that the transition density $Q(u;t|\xi,s)$ for Eq.\eqref{eq:1} meets the following Fokker-Planck equation \cite{Duan}
\begin{equation}\label{eq:8}
\begin{split}
\frac{\partial}{\partial t}Q(x;t|\xi,s)&=-\sum^{d}_{i=1}\frac{\partial}{\partial x_i}(f_i(x)Q(x;t|\xi,s))+\frac{1}{2}\sum^{d}_{i,j=1} \frac{\partial^2}{\partial x_i x_j}(\eta^{i,j}(x)Q(x;t|\xi,s)),
\end{split}
\end{equation}
where $\eta(x)=\sigma(x)\sigma(x)^T$.

Finally, we will further explain the significance of the maximum likelihood transition path. In Eq.\eqref{eq:7}, we can see that the density function $\mathcal{P}_{A}(x,t)$ is determined by the transition probability density $Q(u;t|\xi;s)$ of the state at three different moments, and there is a peak of $\mathcal{P}_{A}(x,t)$ at $t\in[0,t_f]$ according to condition $A$. However, the stochastic orbit $X(t)$ will reach a position state $x_{m}(t)$ of maximum likelihood when the peak value is reached. Therefore, we can maximize the density function $\mathcal{P}_{A}(x,t)$ to get the state $x_m(t)$,
\begin{equation}\label{eq:9}
x_{m}(t)=\arg\max_{x}\mathcal{P}_{A}(x, t).
\end{equation}
We can obtain the state $x_{m}(t)$ through global numerical optimization method. In addition, it is worth mentioning that the most probable state depends on the choice of $t_f$, that is, different values of $t_f$, the maximum likelihood transition path will change.

\section{Numerical experiments}
A few examples verifying the feasibility of our method for detecting the maximum likelihood transition path from data will be presented. Here, we consider two stochastic dynamic systems driven by additive and multiplicative Gaussian noise, respectively. Basic functions are selected as polynomials.

\textbf{Example 1}: The double-well system driven by additive Gaussian noise
\begin{equation}\label{eq:10}
%\begin{aligned}
dX_{t}=(4X_t-X_t^3)dt+{dB_{t}},
%\end{aligned}
\end{equation}
where $X_t$ is a $\mathbb{R}$-valued stochastic process and $\sigma=1$. Firstly, we select the basis function $\Theta(X)$, which has different sizes $n$, and then obtain the cross-validation score $\delta_{\Theta}$ by calculating Eq.\eqref{eq:5}. Figure $1$ shows the relationship between $n$ and the cross validation score: as $n$ increases, the cross validation decreases gradually, and tends to constant.
\begin{figure}[!h]
\subfigure{
\begin{minipage}[t]{0.5\linewidth}
%\centering{$\epsilon=0.01$}
\centerline{\includegraphics[height = 6cm, width = 6cm]{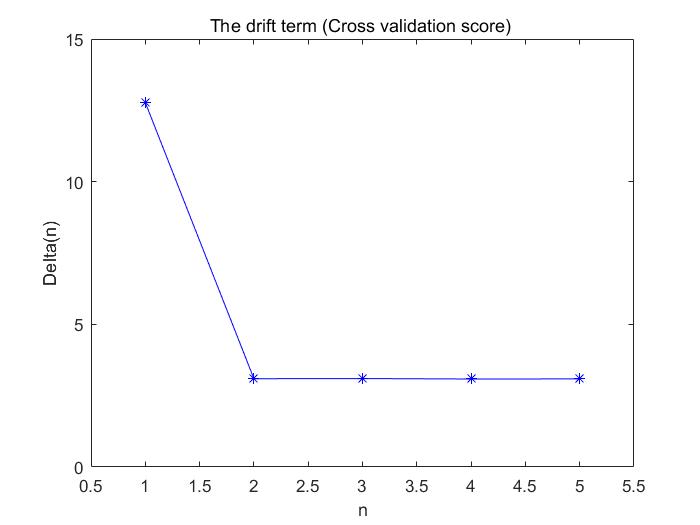}}
%\caption{fig1}
\end{minipage}%
}%
\subfigure{
\begin{minipage}[t]{0.5\linewidth}
%\centering{$\epsilon=0.1$}
\centerline{\includegraphics[height = 6cm, width = 6cm]{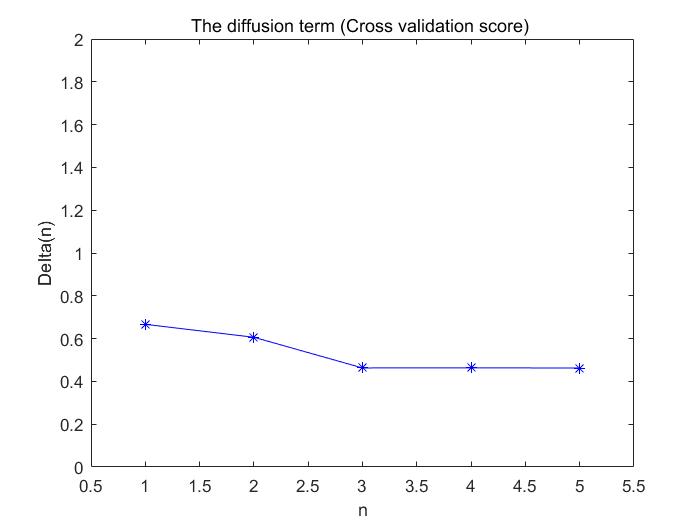}}
%\caption{fig2}
\end{minipage}%
}%
%\centering
\caption{The relationship between size $n$ and the cross validation score.}
\label{fig:1}
\end{figure}

According to the cross validation score,  we can determine that the size of the dictionary $\Theta(X)$ of the drift term and the diffusion term are $n=2$ and $n=1$, respectively. And the results for the Eq.\eqref{eq:4} are as follows,\\
\begin{minipage}{\textwidth}
 \begin{minipage}[t]{0.45\textwidth}
  \centering
     \makeatletter\def\@captype{table}\makeatother\caption{The drift}
       \begin{tabular}{cccc}
         \hline
         basis & $f(x)$ \\
           1   &    0\\
          $x$  &    4.0362   \\
         $x^2$ &    0 \\
         $x^3$ &    -1.0225 \\
         $x^4$ &    0\\
         $x^5$ &    0\\
         \hline
	\end{tabular}
  \end{minipage}
  \begin{minipage}[t]{0.45\textwidth}
   \centering
        \makeatletter\def\@captype{table}\makeatother\caption{The diffusion}
         \begin{tabular}{cccc}
           \hline
           basis & $\sigma^2$ \\
             1   &    1.0094\\
            $x$  &    0 \\
           $x^2$ &    0 \\
           $x^3$ &    0 \\
           $x^4$ &    0 \\
           $x^5$ &    0 \\
           \hline
	  \end{tabular}
   \end{minipage}
\end{minipage}
So, the results corresponding to the original equation \eqref{eq:10} are $f(x)=4.0362x-1.0225x^3$ and $\sigma^2=1.0094$.

Next, we study the dynamic behavior of the double-well system with additive Gaussian noise and draw the maximum likelihood transition path diagram. In Figure $2$, (a) gives the maximum likelihood transition path diagram of original equation under what conditions, and (b) gives the maximum likelihood transition path of double-well system obtained from learning data under the same conditions.

\begin{figure}[!h]
\subfigure{
\begin{minipage}[t]{0.5\linewidth}
%\centering{$\epsilon=0.01$}
\centerline{\includegraphics[height = 6cm, width = 6cm]{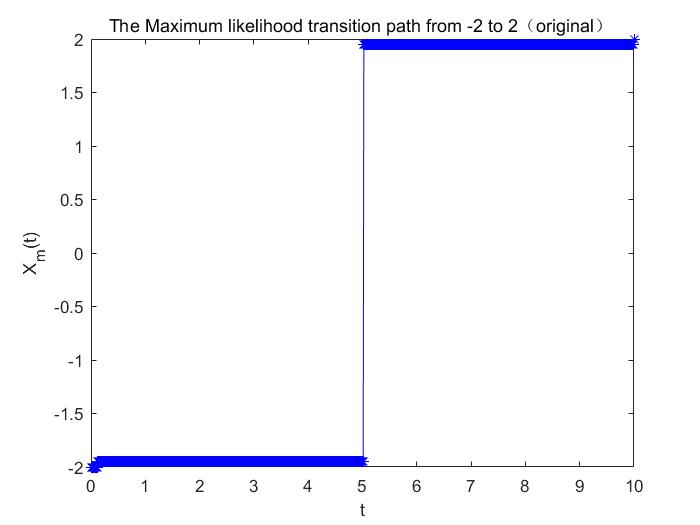}}
%\caption{fig1}
\end{minipage}%
}%
\subfigure{
\begin{minipage}[t]{0.5\linewidth}
%\centering{$\epsilon=0.1$}
\centerline{\includegraphics[height = 6cm, width = 6cm]{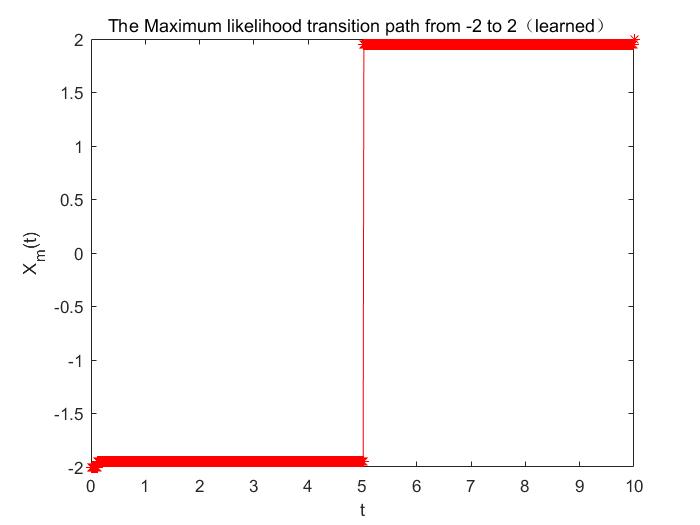}}
%\caption{fig2}
\end{minipage}%
}%
%\centering
\caption{Learned the maximum likelihood transition path for learned double-well system \eqref{eq:10} is on the right, and the original system is on the left.}
\label{fig:2}
\end{figure}

\textbf{Example 2}: The double-well system driven by multiplicative Gaussian noise
\begin{equation}\label{eq:11}
%\begin{aligned}
dX_{t}=(4X_t-X_t^3)dt+(X_t+1){dB_{t}},
%\end{aligned}
\end{equation}
where $X_t$ is a $\mathbb{R}$-valued stochastic process and $\sigma(x)=x+1$. Firstly, we select the basis function $\Theta(X)$, which has different sizes $n$, and then obtain the cross-validation score $\delta_{\Theta}$ by calculating Eq.\eqref{eq:5}. Figure $3$ shows the relationship between $n$ and the cross validation score: as $n$ increases, the cross validation decreases gradually, and tends to constant.
\begin{figure}[!h]
\subfigure{
\begin{minipage}[t]{0.5\linewidth}
%\centering{$\epsilon=0.01$}
\centerline{\includegraphics[height = 6cm, width = 6cm]{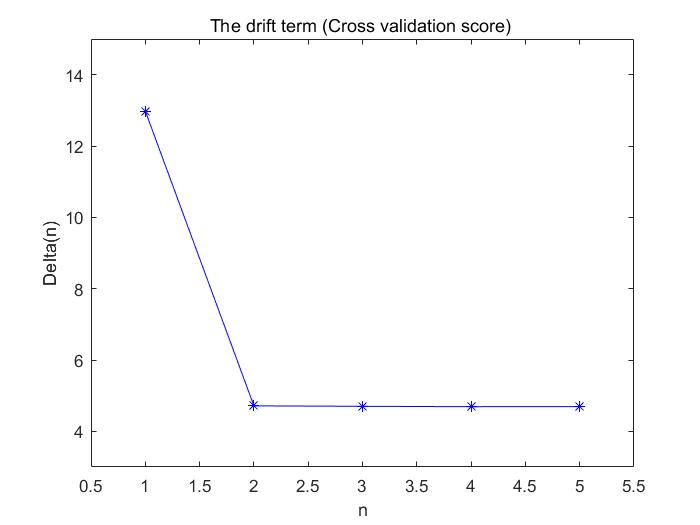}}
%\caption{fig1}
\end{minipage}%
}%
\subfigure{
\begin{minipage}[t]{0.5\linewidth}
%\centering{$\epsilon=0.1$}
\centerline{\includegraphics[height = 6cm, width = 6cm]{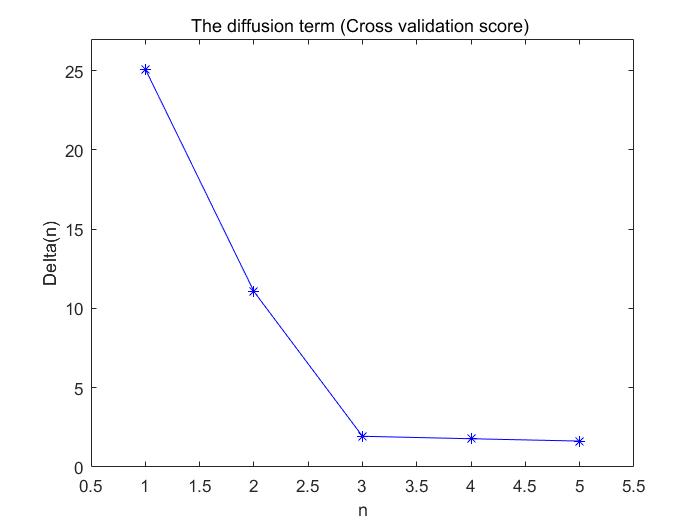}}
%\caption{fig2}
\end{minipage}%
}%
%\centering
\caption{The relationship between size $n$ and the cross validation score.}
\label{fig:3}
\end{figure}

According to the cross validation score, we can determine that the size of the dictionary $\Theta(X)$ of the drift term and the diffusion term are $n=2$ and $n=1$, respectively. And the results for the Eq.\eqref{eq:4} are as follows,\\
\begin{minipage}{\textwidth}
 \begin{minipage}[t]{0.45\textwidth}
  \centering
     \makeatletter\def\@captype{table}\makeatother\caption{The drift}
       \begin{tabular}{cccc}
         \hline
         basis & $f(x)$ \\
           1   &    0\\
          $x$  &    3.9464   \\
         $x^2$ &    0 \\
         $x^3$ &    -0.9998 \\
         $x^4$ &    0\\
         $x^5$ &    0\\
         \hline
	\end{tabular}
  \end{minipage}
  \begin{minipage}[t]{0.45\textwidth}
   \centering
        \makeatletter\def\@captype{table}\makeatother\caption{The diffusion}
         \begin{tabular}{cccc}
           \hline
           basis & $\sigma^2(x)$ \\
             1   &    1.1450\\
            $x$  &    1.9535 \\
           $x^2$ &    0.9135 \\
           $x^3$ &    0 \\
           $x^4$ &    0 \\
           $x^5$ &    0 \\
           \hline
	  \end{tabular}
   \end{minipage}
\end{minipage}
So, the results corresponding to the original equation \eqref{eq:11} are $f(x)=3.9464x-0.9998x^3$ and $\sigma^2(x)=1.1450+1.9535x+0.9135x^2$.

Next, we study the dynamic behavior of the double-well model with multiplicative Gaussian noise and draw the maximum likelihood transition path diagram. In Figure $4$, (a) gives the maximum likelihood transition path diagram of original equation under what conditions, and (b) gives the maximum likelihood transition path of double-well equation obtained from learning data under the same conditions.

\begin{figure}[!h]
\subfigure{
\begin{minipage}[t]{0.5\linewidth}
%\centering{$\epsilon=0.01$}
\centerline{\includegraphics[height = 6cm, width = 6cm]{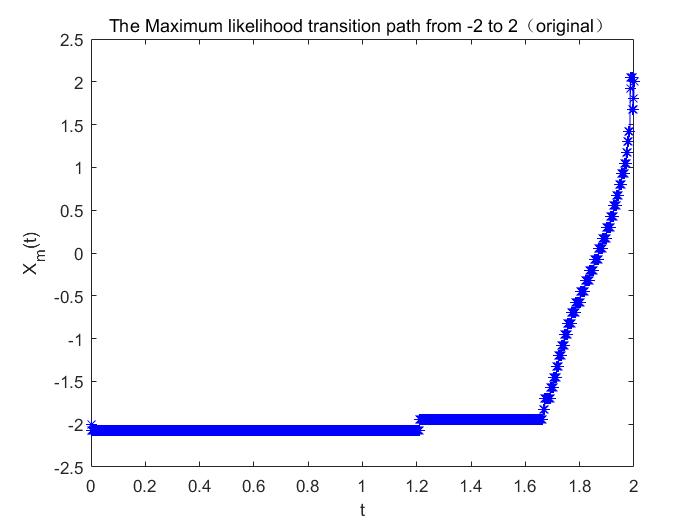}}
%\caption{fig1}
\end{minipage}%
}%
\subfigure{
\begin{minipage}[t]{0.5\linewidth}
%\centering{$\epsilon=0.1$}
\centerline{\includegraphics[height = 6cm, width = 6cm]{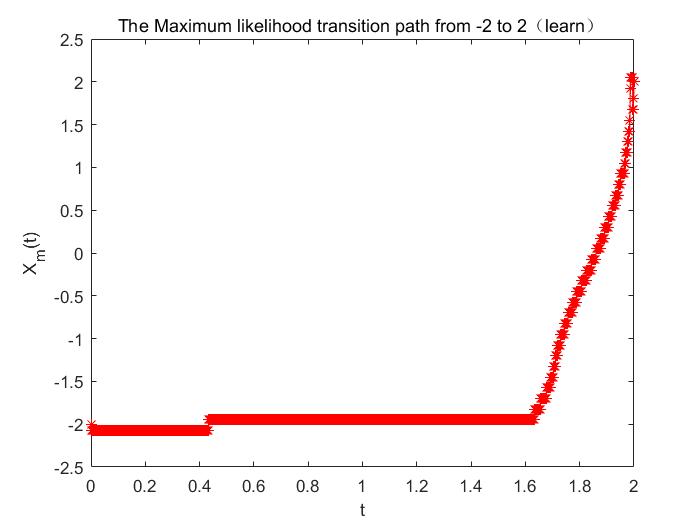}}
%\caption{fig2}
\end{minipage}%
}%
%\centering
\caption{%\textbf{The effect of  $\alpha$-stable L\'evy noise  intensities on  the maximum likelihood transition path for double-well system} .
Learned the maximum likelihood transition path for learned double-well system \eqref{eq:11} is on the right, and the original system is on the left.}
\label{fig:4}
\end{figure}

\section{Conclusion and discussions}
In this work, we proposed a way to obtain dynamical quantities of stochastic dynamical systems from data. Especially, we demonstrated how to extract the maximum likelihood transition path from data. Another contribution is that our method adopted the Kramers-Moyal formula and extended SINDy to deal with stochastic differencial equations with multiplicative Gaussian noise. The advantage of our method is that it might be used to analyze complex phenomena directly from data and obtain quantitative characterization of stochastic systems, which
reduces our dependence on models. Furthermore, our method has the potential to be used to analyze actual data such as genetic data, climate data, etc.

In addition, there are some interesting extensions for this study. First, we only consider the one-dimensional case. For the high-dimensional case, the coordinate binning method is not good for removing noise. How to deal with the high-dimensional stochastic dynamic system situation will become our challenge. On the other hand, for non-Gaussian noise, such as L\'evy noise, it is also a subject worthy of our consideration. Currently, we are conducting research on L\'evy noise.

\section*{Acknowledgements}
We would like to thank Prof. Xiangjun Wang, Dr. Xiujun Cheng, Dr. Xiaoli Chen, Dr. Wei Wei, Dr. Zibo Wang and Dr. Li Lv for helpful discussions. This work was supported by the NSFC grants 11531006, 11801192, 11771449 and NSF grant 1620449.

\section*{Appendix}
\label{SI}
\renewcommand{\theequation}{SM.\arabic{equation}}
\renewcommand{\thefigure}{S.\arabic{figure}}
\setcounter{equation}{0}
\setcounter{figure}{0}

\medskip
\noindent \textbf{S1.  The Kramers-Moyal expansion}

We know that the solution of SDE \eqref{eq:1} is an It$\mathrm{\hat{o}}$ diffusion, which has an identity of the form
\begin{equation}\label{S1:1}
\begin{split}
p(x,t+\Delta{t})=\int{p(x^{'},t)p(x,t+\Delta{t}|x^{'},t)dx^{'}}
\end{split}
\end{equation}
with the probability density $p(x,t+\Delta{t})$ and $p(x^{'},t)$ at time $t+\Delta{t}$ and $t$, respectively, and transition density $p(x,t+\Delta{t}|x^{'},t)$. Firstly, suppose that all the conditional moments $M^{(n)}(x^{'},t,\Delta{t})$ exist and have the following form
\begin{equation}\label{S1:2}
\begin{split}
M^{(n)}(x^{'},t,\Delta{t})&=\mathbb{E}{[{(x(t+\Delta{t})-x(t))^{n}}|x(t)=x]}\\
&=\int{(x-x^{'})^{n}{p(x,t+\Delta{t}|x^{'},t)dx}}.
\end{split}
\end{equation}
In addition, the transition density $p(x,t+\Delta{t}|x^{'},t)$ can be written as
\begin{equation}\label{S1:3}
\begin{split}
p(x,t+\Delta{t}|x^{'},t)&=\int{\delta(y-x)p(y,t+\Delta{t}|x^{'},t)dy}\\
&=\sum_{n=0}^{\infty}{\frac{1}{n!}(-\frac{\partial}{\partial{x}})^{n}}\int{(y-x^{'})^{n}{p(y,t+\Delta{t}|x^{'},t)dy}\delta(x^{'}-x)}\\
&=\Bigg[1+\sum_{n=1}^{\infty}{\frac{1}{n!}(-\frac{\partial}{\partial{x}})^{n}}M^{(n)}(x^{'},t,\Delta{t})\Bigg]\delta(x^{'}-x),
\end{split}
\end{equation}
where
\begin{equation*}
\begin{split}
\delta(y-x)&=\delta(y-x^{'}+x^{'}-x)\\
&=\sum_{n=0}^{\infty}{\frac{(y-x^{'})^{n}}{n!}(-\frac{\partial}{\partial{x^{'}}})^{n}}\delta(x^{'}-x)\\
&=\sum_{n=0}^{\infty}{\frac{(y-x^{'})^{n}}{n!}(-\frac{\partial}{\partial{x}})^{n}}\delta(x^{'}-x).
\end{split}
\end{equation*}
Substituting \eqref{S1:3} into \eqref{S1:1}, we then have
\begin{equation}\label{S1:4}
\begin{split}
p(x,t+\Delta{t})-p(x,t)&=\frac{\partial{p(x,t)}}{\partial{t}}\Delta{t}+o(\Delta{t}^2)\\
&=\int{p(x,t+\Delta{t}|x^{'},t)p(x^{'},t)dx^{'}}-p(x,t)\\
&=\int\Bigg[1+\sum_{n=1}^{\infty}{\frac{1}{n!}\Bigg(-\frac{\partial}{\partial{x}}\Bigg)^{n}}M^{(n)}(x^{'},t,\Delta{t})\Bigg]\delta(x^{'}-x)p(x^{'},t)dx^{'}-p(x,t)\\
%&=\sum_{n=1}^{\infty}{\frac{1}{n!}\Bigg(-\frac{\partial}{\partial{x}}\Bigg)^{n}}M^{(n)}(x^{'},t,\tau)\delta(x^{'}-x)p(x^{'},t)dx^{'}\\
&=\sum_{n=1}^{\infty}\Bigg(-\frac{\partial}{\partial{x}}\Bigg)^{n}\Bigg[\frac{M^{(n)}(x^{'},t,\Delta{t})}{n!}\Bigg]p(x,t).
\end{split}
\end{equation}
Assume that the Taylor series of $M^{(n)}(x^{'},t,\Delta{t})$ is $n!(T^{(n)}(x,t)\Delta{t}+o(\Delta{t}^2))$. Next, dividing both sides of formula \eqref{S1:4} by $\Delta{t}$ at the same time and let $\Delta{t}\rightarrow{0}$,
\begin{equation}\label{S1:5}
\begin{split}
\frac{\partial{p(x,t)}}{\partial{t}}=\sum_{n=1}^{\infty}\Bigg(-\frac{\partial}{\partial{x}}\Bigg)^{n}[T^{(n)}(x,t)p(x,t)].%=\mathcal{L}_{KM}p(x,t)
\end{split}
\end{equation}
%with $\mathcal{L}_{KM}=\sum_{n=1}^{\infty}(-\frac{\partial}{\partial{x}})^{n}D^{(n)}(x,t)$.
So far, we have got the Kramers-Moyal expansion and then the integral notation of the SDE \eqref{eq:1} is
\begin{equation}\label{S1:6}
\begin{split}
X(t)=x+\int_0^{t}f(X(s))ds+\int_0^t\sigma(X(s))dW_s.
\end{split}
\end{equation}
Suppose that $f,\sigma$ have the expansion
\begin{equation}\label{S1:7}
\begin{split}
&f(X(t))=f(x)+f^{'}(x)(X(t)-x)+\cdots\\
&\sigma(X(t))=\sigma(x)+\sigma^{'}(x)(X(t)-x)+\cdots.
\end{split}
\end{equation}
Substituting \eqref{S1:7} into \eqref{S1:6},
\begin{equation}\label{S1:8}
\begin{split}
X(t)-x=&\int_0^{t}{f(x)}ds+\int_0^{t}f^{'}(x)(X(s)-x)ds+\cdots\\
&+\int_0^{t}{\sigma(x){dW_s}}+\int_0^{t}{\sigma^{'}(x)(X(s)-x){dW_s}}+\cdots\\
=&\int_0^{t}{f(x)}ds+\int_0^{t}f^{'}(x)\int_0^s{f(x_{u})d{u}}ds+\int_0^{t}f^{'}(x)\int_0^s{\sigma(x_{u})dW_{u}}ds+\cdots\\
&+\int_0^{t}{\sigma(x){dW_s}}+\int_0^{t}{\sigma^{'}(x)\int_0^s{f(x_{u})d{u}}{dW_s}}+\int_0^{t}{\sigma^{'}(x)\int_0^s{\sigma(x_{u})dW_{u}}{dW_s}}+\cdots,
\end{split}
\end{equation}
we then obtain
\begin{equation}\label{S1:9}
\begin{split}
\mathbb{E}[X(t)-x]&=tf(x)+\frac{t^2}{2}f^{'}(x)f(x)+\mathbb{E}\Bigg[\int_0^{t}{\sigma^{'}(x)\int_0^s{\sigma(x_{u})dW_{u}}dW_s}\Bigg]\\
&=tf(x)+\frac{t^2}{2}f^{'}(x)f(x)+\sigma^{'}(x)\sigma(x)\mathbb{E}\Bigg[\int_0^{t}{W_sdW_s}\Bigg].
\end{split}
\end{equation}
In the sense of It$\mathrm{\hat{o}}$, the expectation $\mathbb{E}[\int_0^{t}{W_sdW_s}]=0$. Therefore,
\begin{equation}\label{S1:10}
\begin{split}
T^{(1)}=\lim_{t\rightarrow{0}}\frac{1}{t}\mathbb{E}[X(t)-x]=f(x).
\end{split}
\end{equation}
In addition,
\begin{equation}\label{S1:11}
\begin{split}
T^{(2)}=&\frac{1}{2}\lim_{t\rightarrow{0}}\frac{1}{t}{\mathbb{E}[(X(t)-x)^2]}\\
=&\frac{1}{2}\lim_{t\rightarrow{0}}\frac{1}{t}\mathbb{E}\Bigg[{\Bigg(\int_0^{t}{f(x)}ds+\int_0^{t}f^{'}(x)\int_0^s{f(x_{u})d{u}}ds+\int_0^{t}f^{'}(x)\int_0^s{\sigma(x_{u})dW_{u}}ds+\cdots}\\
&{+\int_0^{t}{\sigma(x){dW_s}}+\int_0^{t}{\sigma^{'}(x)\int_0^s{f(x_{u})d{u}}{dW_s}}+\int_0^{t}{\sigma^{'}(x)\int_0^s{\sigma(x_{u})dW_{u}}{dW_s}}+\cdots}\Bigg)^2\Bigg]\\
=&\frac{1}{2}\lim_{t\rightarrow{0}}\frac{1}{t}\Bigg(t^2f^2(x)+\frac{t^4}{4}(f^{'}(x))^{2}f^2(x)+(f^{'}(x))^{2}\sigma^2(x)\mathbb{E}\Bigg(\int_0^tW_sds\Bigg)^2+\sigma^2(x)\mathbb{E}(W_t^2)\\
&+t^2(\sigma^{'}(x))^2f^2(x)\mathbb{E}(W_t^2)+(\sigma^{'}(x))^2\sigma^2(x)\mathbb{E}\Bigg(\int_0^tW_sdW_s\Bigg)^2+\cdots\\
&+\frac{t^3}{2}f^{2}(x)f^{'}(x)+tf(x)f^{'}(x)\sigma(x)\mathbb{E}\Bigg(\int_0^tW_sds\Bigg)+f^{'}(x)\sigma^2(x)\mathbb{E}\Bigg(W_t\int_0^tW_sds\Bigg)\\
&+tf^{'}(x)\sigma(x)\sigma^{'}(x)f(x)\mathbb{E}\Bigg(W_t\int_0^tW_sds\Bigg)+t\sigma(x)\sigma^{'}(x)f(x)\mathbb{E}(W_t^2)+\cdots\Bigg)\\
=&\frac{1}{2}\sigma^{2}(x).
\end{split}
\end{equation}
Similarly, $T^{(n)}=\frac{1}{n!}\lim_{t\rightarrow{0}}\frac{1}{t}{\mathbb{E}[(X(t)-x)^n]}=0,~n\geq{3}$. So, the Kramers-Moyal formula is obtained.

\medskip
\noindent \textbf{S2.  The conditional density function}

Suppose that SDE \eqref{eq:1} has a unique strong solution that has a strictly positive probability density, and the conditional density function of the solution also exists. In fact, we know that the value of the conditional density function is only related to the immediately preceding moment through the Markov property of SDE \eqref{eq:1}, therefore, for $0<t<t_f$, we have
\begin{equation}\label{p8}
\begin{split}
p(X(t_f)=x_{f}|X(0)=x_0;X(t)=x)&=p(X(t_f)=x_{f}|X(t)=x)\\
&=Q(x_{f},t_f|x,t).
\end{split}
\end{equation}
In addition, if conditions $X(0)=x_0$ and $X(t_f)=x_{f}$ are given, according to Bayesian formula, we can calculate the conditional density function of the solution
\begin{equation}\label{p9}
\begin{split}
p(X(t)&=x|X(0)=x_0;X(t_f)=x_{f})\\
&= \frac{p(X(t)=x;X(0)=x_0; X(t_f)=x_{f})}{p(X(0)=x_0;X(t_f)=x_{f})}\\
&= \frac{p(X(t)=x|X(0)=x_0)p(X(t_f)=x_{f}|X(t)=x)}{p(X(t_f)=x_{f}|X(0)=x_0)}.
\end{split}
\end{equation}
We then substitute  \ref{p8} into \ref{p9}, and the formal expression of conditional density function $\mathcal{P}_{A}(x, t)$ will be obtained
\begin{equation}\label{p10}
\begin{split}
\mathcal{P}_{A}(x,t)&=p(X(t)=x|X(0)=x_0;X(t_f)=x_{f})\\
&=\frac{Q(x_{f},t_f|x,t)Q(x,t|x_0,0)}{Q(x_{f},t_f|x_0,0)} .
\end{split}
\end{equation}
%Thus, the conditional probability density is obtained.

\end{document}